\documentclass[nameyear,seceqn,dvips]{arx_bjps}
\usepackage{graphicx}

\volume{24}
\issue{3}
\pubyear{2010}
\firstpage{479}
\lastpage{501}
\doi{10.1214/09-BJPS105}

\makeatletter
\DeclareMathAlphabet\mathcaligr{OMS}{cmsy}{m}{n}
\newtheorem{theorem}{Theorem}[section]
\newtheorem{lemma}{Lemma}[section]

\makeatother

\begin{document}
\begin{frontmatter}

\title{Reduced long-range dependence combining Poisson bursts with
on--off sources}
\runtitle{Reduced LRD combining Poisson bursts with on--off sources}

\begin{aug}
\author{\fnms{David A.}\inits{D. A.}
\snm{Rolls}\ead[label=e1]{drolls@unimelb.edu.au}} 

\runauthor{D. A. Rolls}
\affiliation{University of Melbourne}

\address{Department of Mathematics and Statistics\\
University of Melbourne\\
Parkville, VIC 3010\\ Australia\\\printead{e1}} 

\end{aug}

\received{\smonth{9} \syear{2008}}
\accepted{\smonth{5} \syear{2009}}

%
\begin{abstract}
A workload model using the infinite source Poisson model for
bursts is combined with the on--off model for within burst activity.
Burst durations and on--off durations are assumed to have heavy-tailed
distributions with infinite variance and finite mean. Since the number
of bursts is random, one can consider limiting results based on
``random centering'' of a random sum for the total workload from all sources.
Convergence results are shown to depend on the tail indices of both the
on--off durations and the lifetimes distributions. Moreover, the
results can be separated into cases depending on those tail indices. In
one case where all distributions are heavy tailed it is shown that the
limiting result is Brownian motion. In another case, convergence to
fractional Brownian motion is shown, where the Hurst parameter depends
on the heavy-tail indices of the distribution of the on, off and burst
durations.
\end{abstract}

%
\begin{keyword}
\kwd{Teletraffic}
\kwd{fractional Brownian motion}
\kwd{on--off process}
\kwd{heavy tails}
\kwd{long range dependence}.
\vspace*{-9pt}
\end{keyword}

\end{frontmatter}
%

\section{Introduction}\label{sec:intro}

Workload models for packet traffic have been described previously
[\citeauthor{brichet96} (\citeyear{brichet96,brichet00});
\citet{levy00}; \citet{mandelbrot69}; \citet{tl86}; \citet{taqqu97}; \citet{willinger97}; \citet{kurtz96}].
An important model has been the strictly alternating on--off model with
heavy tailed on or off times [\citet{taqqu97}; \citet{willinger97}], in which
each source creates work at constant rate for all time. Another
important model has been the ``$M/G/\infty$ queue'' described by
\citet{cox}, also called the ``infinite source Poisson'' model [\citet
{mikosch2002}]. In this model, traffic arrives as independent bursts of
heavy-tailed size or duration at the time points of a homogeneous
Poisson process. There is no variability within each burst and bursts
are considered to have the same constant rate. (See \citet
{mikosch2002} for a treatment of convergence results in both models.)
The use of heavy tails is motivated by, for example, empirical evidence
on the sizes of WWW objects [\citet{crovella98}].
For these models, the cumulative work centered about its mean is shown
to be a fractional Brownian motion (fBm) in appropriate limiting
regimes, although other limiting regimes have been studied which show
convergence to stable L\'{e}vy motion [e.g., \citet{taqqu97}; \citet{mikosch2002}; \citet{kaj07}].

A number of more recent models combine Poisson arrivals of bursts with
assumptions on the burst volume, burst rate, or by introducing some
dynamics within the burst. In \citet{maulik02} the infinite source
Poisson model is used but with an independent random rate given to each
burst. There is no variability within bursts. In \citet{maulik03}
bursts are given a ``transmission schedule'' according to an
$H$-self-similar process with nondecreasing, cadlag paths (necessarily
$H \geq1$), and a ``volume'' of data (work) whose distribution has
infinite variance and finite mean. Under a ``fast growth'' condition
they show convergence (of the finite-dimensional distributions) of the
cumulative work, centered by its mean, to a Gaussian self-similar
process that generally lacks stationary increments. So, their limit is
generally not fractional Gaussian noise. In \citet{caglar04} ``flows''
arrive according to a Poisson process, have infinite variance
Pareto-distributed holding times, but packets within a flow arrive
according to a compound Poisson process, and packet sizes have finite
variance. The increment of the cumulative workload is shown to be
fractional Gaussian noise, with Hurst parameter $H$ depending on
holding time.

Cluster models have been proposed to model the packet arrival process.
While not modeling the workload per se (absent a model for packet
size), they also model variability within a ``flow.'' In \citet{hohn03}
and \citet{fay06} flows arrive according to a Poisson process. The
number of packets in each flow and the times between packets are random
and may have infinite variance. Long-range dependence arises when the
number of packets has infinite variance and finite mean.

This paper studies the workload generated in a model where sources
arrive at Poisson time points and have heavy-tailed ``session''
durations. During its session a source is stationary with independent
on and off durations. The on times are independent and identically
distributed (i.i.d.), as are the off times, and at least one of the two
distributions is assumed to be heavy tailed with infinite variance but
finite mean. Here ``session'' is a euphemism for a structure above that
of the on--off behavior. Thus the model here can be considered as a
hybrid of the infinite source Poisson model and the alternating on--off
model, although the proof here owes more to the latter. The model
studied here simplifies the one studied by \citet{rolls2003} without
changing the main results, by assuming the number of sessions has
achieved stationarity.

The approach used in this paper can establish a (traditional) workload
limit result where the variability of the workload about its mean is
fractional Brownian motion under certain assumptions. Moreover, the
Hurst parameter is the same as that from a corresponding infinite
source Poisson model (i.e., from the model if there was no variability
within bursts). But the main point here is something else. The on--off
behavior imposes additional variability on top of that from the
sessions themselves. Under so-called ``random centering'' this
additional variability is shown to have a corresponding Hurst parameter
that depends on the indices of both the heavy-tailed session durations
and the on--off durations. The formula for the Hurst parameter is a new
expression, different from the standard limiting results from either
the on--off or infinite source Poisson models. This might arise if
session start and end times were announced and this information was
incorporated into the centering. Two cases are considered, giving rise
to either a fractional Brownian motion or a Brownian motion. It is
noteworthy than under random centering, the limit can be Brownian
motion even when both the session durations and the on or off durations
are heavy tailed.

The rest of this paper is organized as follows. Section \ref
{background:sec} provides background and introduces necessary notation.
Section \ref{modeldesc:sec} describes the alternating on--off model
with session lifetimes. Section \ref{randomcentering:sec} establishes
convergence results for this model. Section \ref{weakconv:sec} extends
the convergence results to weak convergence. Section~\ref
{simresults:sec} contains the result from several simulations. Section
\ref{conc:sec} provides conclusions and discusses possible future work.

\section{Background}
\label{background:sec}

For the on--off processes we use the notation of \citet{taqqu97}. Let
$\{W(t),  t\geq0\}$ be a strictly alternating, stationary on--off
process that takes the value 0 during an off period and 1 during an on
period. Let $\{W_i(t),  t\geq0\}$, $i=1, 2, \ldots$ be mutually
independent copies of $\{W(t),  t\geq0\}$. Define the autocovariance
and mean of $W(t)$ by
\[
r(t) =\mathbb{E}[W^2(t)] - (\mathbb{E}[W(t)])^2\quad  \mbox{and}\quad \mu
_W=\mathbb{E}[W(t)],
\]
respectively.

The on and off times are nonnegative and independent of each other.
Also, the lengths of the on periods are i.i.d. as are the off periods.
Let $f_j(x),  F_j(x)$, and~$\mu_j,$ be the density, distribution
function, and mean for the duration of on periods (\mbox{$j=1$}) or off
periods ($j=2$), respectively.

The on and off times are assumed to be heavy tailed with regularly
varying tails: $\bar{F}_{j}(x)=1-F_j(x)\sim x^{-\alpha_j}L_j(x)$ as
$x \rightarrow\infty$, $1<\alpha_j<2$,
where $L_j>0$ is a slowly varying function, $j=1,2$. (In fact, one can
assume only one distribution is heavy tailed. The results would be
unchanged.) Then
\[
\mu_W= \frac{\mu_1}{\mu_1+\mu_2} \quad  \mbox{and}\quad \sigma_W^2=\operatorname{Var}\bigl[ W(t)-\mathbb{E}[W(t)]\bigr]=\frac{\mu_1\mu_2}{(\mu_1
+ \mu_2)^2}.
\]
Finally, set $\alpha_\mathit{{min}}=\min\{\alpha_1,\alpha_2\}$ and to
specify the indices, let $(\mathit{min},  \mathit{max})= (1,2)$ if $\alpha_1<\alpha_2$
and $(\mathit{min},  \mathit{max})= (2,1)$ if $\alpha_2<\alpha_1$. Under these
assumptions Taqqu, Willinger and Sherman (\citeyear{taqqu97}) showed
\[
\mathcaligr{L}\lim_{T\rightarrow\infty} \mathcaligr{L} \lim
_{M\rightarrow\infty} \frac{(\int_0^{Tt}( \sum_{i=1}^M
W_i(u)) \, du -TM\mu_Wt)}{T^H L^{1/2}(T)M^{1/2}} = \sigma
_{\mathit{lim}}Z_H(t),
\]
where $\mathcaligr{L}\lim$ means convergence of the finite-dimensional
distributions, $Z_H(t)$ is standard fractional Brownian motion,
$H=(3-\alpha_\mathit{{min}})/2$ and
%
\begin{equation}
\label{sigmalim1}
\sigma_{\mathit{lim}}^2=\frac{2 \mu_{\mathit{max}}^2}{\mu_W^3(\alpha_\mathit{{min}}-1)
(3-\alpha_\mathit{{min}})(2-\alpha_\mathit{{min}})}.
\end{equation}
(See \citet{taqqu97} for $\alpha_1=\alpha_2$.)

Conditions so that $r(t)$ has a tail regularly varying at infinity are
given by \citet{brichet00} and \citet{heath98}. For the results here
it is assumed that
%
\begin{equation}
\label{regvarfuncr}
r(u) = \frac{\sigma_{\mathit{lim}}^2}{2} (3-\alpha_\mathit{{min}}) (2-\alpha_\mathit{{min}})
u^{1-\alpha_\mathit{{min}}}L_r(u)
\end{equation}
for some slowly varying function $L_r(u)$.

Analogous convergence results to fractional Brownian motion have been
obtained for the infinite source Poisson model [\citet{kurtz96}]. In
those results it is not the number of sources, $M$, but rather the
arrival rate of a homogeneous Poisson process, say $\lambda$, that
goes to infinity, followed by the time rescaling $T$.

\section{The model description}
\label{modeldesc:sec}

In this section we define the alternating on--off model with
heavy-tailed session lifetimes. For the on--off processes we use the
same notation and assumptions as in Section \ref{background:sec}. In
particular, at least one of the on or off-period distributions is
assumed heavy tailed so that $1< \alpha_\mathit{{min}} <2$. Also, let $\{G(t),
t \geq0 \}$ be the mean-zero Gaussian process with autocovariance
$r(t)$ in \citet{taqqu97} arising as
\[
\mathcaligr{L}\lim_{n \rightarrow\infty} \frac{\sum_{i=1}^n(
W_i(t)-\mu_W ) }{\sqrt{n}} = G(t).
\]

For the sessions, we make assumptions on their arrivals and durations.
Let $\{T_{\lambda_n,i},-\infty<i<\infty\}$ be the arrival times of a
rate $\lambda_n$ Poisson process on $\mathbb{R}$, labeled so that
$T_{\lambda_n,0}<0<T_{\lambda_n,1}$, where $\{\lambda_n,
n=1,2,\ldots\}$ is a sequence of positive constants such that $\lambda
_n \rightarrow\infty$ as $n \rightarrow\infty$. For the durations
of the sessions, let $\{V_i\}$ be an i.i.d. sequence with continuous
distribution $H$ and finite mean $\mu_V$.
Here $V_i$ will be the lifetime (i.e., holding time, session length)
for source $i$.
It is assumed that the distribution function $H(x)$ of the lifetimes is
Lipschitz continuous (e.g., satisfied if $H$ has a bounded derivative),
and has a regularly varying tail so we may write
\[
\bar{H}(x) = x^{-\alpha_{\mathit{sess}}}L_V(x),\qquad
\alpha_{\mathit{sess}}>0,
\]
where $L_V(x)$ is slowly varying. Consideration of the integrated tail
of $H(x)$, $\bar{H}_I(x) = \int_x^\infty\bar{H}(z)\, dz $ will be
necessary. If $\alpha_{\mathit{sess}}>1$, by Karamata's theorem
[Bingham, Goldie and Teugels (\citeyear{bingham87}), p.~28]
\begin{equation}
\label{Hbartail}
\bar{H}_I(x) \sim\frac{x^{-\alpha_{\mathit{sess}}+1}}{\alpha
_{\mathit{sess}}-1}L_V(x)\qquad  \mbox{as } x \rightarrow\infty.
\end{equation}
In the boundary case $\alpha_{\mathit{sess}}=1$, with the additional assumption
that $\int_1^\infty L_V(t)/\break t  \, dt < \infty$ then
\[
\frac{1}{L_V(x)}\int_x^\infty\frac{L_V(t)}{t} \, dt \rightarrow
\infty\qquad  \mbox{as } x\rightarrow\infty.
\]

The processes $\{T_{\lambda_n,i} \}$, $\{V_i \}$ and $\{W_i(t)\}$ are
all assumed independent. Together, the combination of the Poisson
arrival process and the sequence of holding times defines the busy
server process of an $M/G/\infty$ queueing system.

A key idea in this paper is the use of random sums. That is, sums whose
upper limit of summation is random and obeys some convergence result of
its own [\citet{gnedenko96}]. With random sums one must distinguish
between nonrandom centering of the sum, and nonrandom centering of the
summands (really a \textit{random} centering of the sum). For the model
presented here, convergence results for both kinds of centering will
play a role and so both are discussed.

Let
\[
A_n(t)= \int_0^t\sum_{i=-\infty}^{\infty}W_i(u) \mathbf
{1}_{[T_{\lambda_n,i},T_{\lambda_n,i}+V_i)}(u)\, du
\]
and
\[
B_n(t)= \int_0^t\sum_{i=-\infty}^{\infty}\mathbf{1}_{[T_{\lambda
_n,i},T_{\lambda_n,i}+V_i)}(u)\, du,
\]
where $\mathbf{1}_A (x)$ is 1 on the set $A$ and 0 otherwise. With
this notation $A_n(t)$ is the total cumulative work in $[0,t)$ for the
alternating on--off sources with lifetimes. Note that $B_n(t)$ is the
total cumulative work in $[0,t)$ from the infinite source Poisson model
whose ``bursts'' are exactly the ``sessions'' in $A_n(t)$ and sources
are on for the complete burst.

The total cumulative work for the infinite source Poisson model with
nonrandom centering is easy to express
%
\begin{equation}
\label{infinitesourcePoisson}
B_n(Tt)-\mu_V\lambda_n Tt = \int_0^{Tt} \Biggl[ \sum_{i=-\infty
}^{\infty}\mathbf{1}_{[T_{\lambda_n,i},T_{\lambda_n,i}+V_i)}(u) -
\mu_V \lambda_n \Biggr] \, du.
\end{equation}
One can show weak convergence of a rescaled version of this quantity to
fBm with $H=(3-\alpha_{\mathit{sess}})/2$ under suitable assumptions.
Similarly, for the alternating on--off model with lifetimes described
here, with nonrandom centering, one can show weak convergence of a
rescaled version of
\[
A_n(Tt)-\mu_W \mu_V\lambda_n Tt
\]
to fBm with the same Hurst parameter as for (\ref{infinitesourcePoisson}).

The focus of this paper is something different, namely
%
\begin{equation}
\label{randomcentering}
A_n(Tt)-\mu_W B_n(Tt)
= \int_0^{Tt}\sum_{i=-\infty}^{\infty} X_{n,i}(u)  \, du,
\end{equation}
where
%
\begin{equation}
\label{forcelife}
X_{n,i}(t) = [W_i(t)-\mu_W ] \mathbf{1}_{[T_{\lambda
_n,i},T_{\lambda_n,i}+V_i)}(t).
\end{equation}
This quantity captures the variability from the on--off dynamics on top
of the variability (at a larger time scale) from the session arrivals
and departures. The goal now will be to establish the convergence
properties of this process.

\section{Main results}
\label{randomcentering:sec}

\begin{figure}[b]

\includegraphics{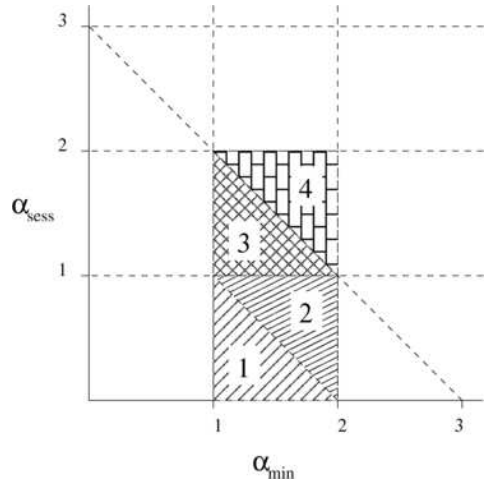}

\caption{Differing cases arising from the
tail indices $\alpha_\mathit{{min}}$ and $\alpha_{\mathit{sess}}$. Case \textup{3} leads to fBm,
while case \textup{4} leads to Brownian motion. The long diagonal line
corresponds to $3-\alpha_{\mathit{min}}-\alpha_{\mathit{sess}}=0$. }\label{lifetimes_cases}
\end{figure}

Here we present our main results, which depend on the tail indices
$\alpha_\mathit{{min}}$ and $\alpha_{\mathit{sess}}$. The results can be separated into
several cases by these values, as shown in Figure~\ref
{lifetimes_cases}. Our first theorem is the limit result for case 3,
while the second theorem is for case 4. In particular, note that
although the tail indices of $\alpha_\mathit{{min}}$ and $\alpha_{\mathit{sess}}$ both
correspond to heavy tails, the limiting result is Brownian motion for
case 4. Cases~1 and 2 are left for future work.

\begin{theorem}
\label{fbmmainthm1}
For $1<\alpha_\mathit{{min}}<2$ and $1<4-\alpha_\mathit{{min}}-\alpha_{\mathit{sess}} <2$
%
\begin{equation}
\mathcaligr{L}\lim_{T \rightarrow\infty}\mathcaligr{L}\lim_{n
\rightarrow\infty} \frac{\int_0^{Tt} \sum_{i=-\infty}^\infty
X_{n,i}(u)  \,du}{\sqrt{\lambda_n L_V(T)L_r(T)}T^H }=
\sigma Z_H(t),
\label{fbmmain1eqn}
\end{equation}
where $H=(4-\alpha_\mathit{{min}}-\alpha_{\mathit{sess}})/2$, $\sigma_{\mathit{lim}}^2$ is
defined by (\ref{sigmalim1}),
\[
\sigma^2= \frac{\sigma_{\mathit{lim}}^2(3-\alpha_\mathit{{min}})(2-\alpha
_\mathit{{min}})}{(\alpha_{\mathit{sess}}-1)(3-\alpha_\mathit{{min}}-\alpha_{\mathit{sess}})(4-\alpha
_\mathit{{min}}-\alpha_{\mathit{sess}})},
\]
and $\{Z_H(t)\}$ is standard fractional Brownian motion with Hurst
parameter $H$.
\end{theorem}

\begin{theorem}
For $1<\alpha_\mathit{{min}}<2$ and $-2<2-\alpha_\mathit{{min}}-\alpha_{\mathit{sess}} <-1,
\alpha_{\mathit{sess}}<2$
\label{fbmmainthm2}
%
\begin{equation}
\label{fbmmainthm2eqn}
\mathcaligr{L}\lim_{T \rightarrow\infty}\mathcaligr{L}\lim_{n
\rightarrow\infty} \frac{\int_0^{Tt} \sum_{i=-\infty} ^\infty
X_{n,i}(u) \, du}{\sqrt{\lambda_n T } }=\sqrt{2c}Z_{0.5}(t),
\end{equation}
where
\[
c= \int_0^\infty r(x) \bar{H}_I(x)\, dx
\]
and $\{Z_{0.5}(t)\}$ is standard Brownian motion.
\end{theorem}

The rest of this section will prove these results in three steps.
Section \ref{incremental.subsec} presents limit results for the
incremental process (i.e, for $n \rightarrow\infty$). Section \ref
{varcovar.subsec} describes the asymptotic behavior of the variance of
the incremental process, which is useful for finding the limiting
process as $T \rightarrow\infty$. Section \ref{mainproof.subsec}
proves the main results.

\subsection{Limit results for the incremental process}
\label{incremental.subsec}

Let $\{Z(t)\}$ be a Gaussian process such that for $s<t$:
%
\begin{eqnarray}\label{myZ}
\mbox{(i)}\hspace*{48pt}  E[Z(t)] &=& 0, \nonumber\\
\mbox{(ii)}\hspace*{40pt}  \operatorname{Var}[Z(t)]&=&\sigma_W^2 \mu_V\quad \mbox{and} \\
\mbox{(iii)}\quad  \operatorname{Cov}[Z(s),Z(t)]&=&r(t-s) \int_{t-s} ^{\infty}
\bar{H}(z)  \,dz, \nonumber
\end{eqnarray}
where $\bar{H}(z)=1-H(z)$. It will be shown that to obtain limiting
results as $T \rightarrow\infty$, $\{Z(t)\}$ takes the place of the
incremental process $\{G(t)\}$ in \citet{taqqu97}.

\begin{theorem}
\label{incrementalprocess}
For the processes $\{X_{n,i}(t)\}$ of (\ref{forcelife}) and $\{Z(t)\}$
of (\ref{myZ})
\[
\mathcaligr{L}\lim_{n \rightarrow\infty} \frac{\sum_{i=-\infty
}^{\infty}X_{n,i}(t)}{\sqrt{\lambda_n}}=Z(t).
\]
\end{theorem}

\begin{pf}
For any $d\in\mathbb{Z}^+$, let $(c_1,\ldots,c_d)$ be arbitrary and
let $0<t_1<t_2<\cdots<t_d$ be a partition of the time axis. To
simplify the exposition, let $I_k=[t_k,t_{k+1})$, $k=1,2,\ldots,d-1$,
$I_0=(-\infty,t_1)$, $I_d=[t_d,\infty)$ and let
\[
Z_n = \frac{c_1 \sum_{k=-\infty}^{\infty} X_{n,k}(t_1)+ \cdots+ c_d
\sum_{k=-\infty}^{\infty}X_{n,k}(t_d)}{\sqrt{\lambda_n}}.
\]

Classify an arrival as type $(i,j)$, $i=0,\ldots,d$, $j=i,\ldots,d$,
if it arrives in $I_i$ and ends in $I_j$. Let $N_{n,i,j}(t)$ be the
number of type $(i,j)$ arrivals by time $t, t>0$. Conditional on an
arrival at time $y$, the probability an arrival is of type $(i,j)$,
$P_{(i,j)}(y)$, is
\[
P_{(i,j)}(y) =\cases{
0, &\quad $ y \notin I_i$, \vspace*{2pt}\cr
H(t_{i+1}-y), &\quad $ y \in I_i$,   $i=j \neq d$, \vspace*{2pt}\cr
H(t_{j+1}-y)-H(t_j-y), &\quad $ y \in I_i$,   $i<j$,   $j \neq d$, \vspace*{2pt}\cr
\bar{H}(t_d-y), &\quad $ y\in I_i$,   $i<j $,  $j=d$, \vspace*{2pt}\cr
1, &\quad $ y \in I_i$,   $i=j=d$.}
\]
By Proposition 5.3 of \citeauthor{ross} (\citeyear{ross}, p.~273),  for fixed $t$, the random
variables $N_{n,i,j}(t)$, $i=0,\ldots,d$, $j=i,\ldots,d$ are mutually
independent and Poisson distributed with known mean (and variance)
%
\begin{eqnarray}
&&\mathbb{E}[N_{n,i,j}(t)] = \operatorname{Var}[N_{n,i,j}(t)] =
\lambda_n \int_0 ^t P_{(i,j)}(y) \,dy = \lambda_n p_{(i,j)} (t), \nonumber \\
\eqntext{i=0,\ldots,d   ,j=i,\ldots,d,}
\end{eqnarray}
where $ p_{(i,j)}(t) = \int_0 ^t P_{(i,j)}(y)\, dy$. (See \citet
{rolls2003} for details.) For each process $k$ of type $(i,j)$ it
contributes to $Z_n$ during its holding period at $t_{i+1},\ldots
,t_j$. Let $(X_{n,k}^{(i,j)}(t_{i+1}), \ldots,
X_{n,k}^{(i,j)}(t_j) ) $ be copies of $ \Big( W(t_{i+1})-\mu
_W, \ldots, W(t_j)-\mu_W \Big) $, mutually independent in $i$, $j$
and $k$. Notice that for fixed $i$ and $j$ the sequence of random
vectors is i.i.d. Then the contribution of the $k$th type $(i,j)$
process is
\[
c_{i+1} X_{n,k} ^{(i,j)}(t_{i+1}) + \cdots+ c_{j} X_{n,k} ^{(i,j)}(t_{j})
\]
and there are $N_{n,i,j}(t_d)$ such processes. Thus we can write
\[
Z_n  \stackrel{D}{=} \sum_{i=0}^{d-1} \sum_{j=i+1} ^d \biggl[\frac
{\sum_{k=1} ^{N_{n,i,j}(t_d)} ( c_{i+1} X_{n,k}
^{(i,j)}(t_{i+1}) + \cdots+ c_{j} X_{n,k} ^{(i,j)}(t_{j}) )}
{\sqrt{\lambda_n}} \biggr].
\]
(Strictly speaking there is also a term from arrivals occurring exactly
at $t_k$, $k=1,2,\ldots,d$. Since these arrivals have probability zero
that term can be safely ignored.)

By Lemmas 8.2.2 and 8.2.3 of \citet{rolls2003} we have
\[
\frac{\sum_{k=1} ^{N_{n,i,j}(t_d)} ( c_{i+1} X_{n,k}
^{(i,j)}(t_{i+1}) + \cdots+ c_{j} X_{n,k} ^{(i,j)}(t_{j}) )
}{\sqrt{\lambda_n}}
\stackrel{D}{\rightarrow} Z^{(i,j)} \qquad  \mbox{as } n \rightarrow\infty,
\]
$ i=0,\ldots,d-1, j=i+1,\ldots,d$, where $Z^{(i,j)}$ are mutually
independent random variables such that
\begin{eqnarray*}
 Z^{(i,j)} &= &\sqrt{p_{(i,j)}(t_d)}\bigl ( c_{i+1} G
^{(i,j)}(t_{i+1}) + \cdots+ c_{j} G ^{(i,j)}(t_{j}) \bigr) \\
& \sim&\mathrm{N}\Biggl(0, p_{(i,j)}(t_d) \Biggl[ \sigma_W^2 c_{i+1}^2
+ \cdots+\sigma_W^2 c_{j}^2 + 2\!\!\sum_{u=i+1}^j \sum_{v=u+1}^j\! c_u
c_v r(|t_v-t_u|) \Biggr]\Biggr),
\end{eqnarray*}
where $r(u)$ is the autocovariance of $\{G(t), t \geq0\}$ and so
\[
Z_n \stackrel{D}{\rightarrow} \sum_{i=0}^{d-1} \sum_{j=i+1} ^d
Z^{(i,j)} \stackrel{D}{=} \sum_{i=1} ^d c_i Z(t_i) \qquad  \mbox{as } n
\rightarrow\infty.
\]
Here $\{G^{(i,j)}(t),t \geq0 \}$ are independent copies of $\{G(t),t
\geq0 \}$ in $i$ and $j$. Since $(c_1,\ldots,c_d)$ and $(t_1,\ldots
,t_d)$ are arbitrary, the result follows from the Cram\'{e}r--Wold theorem.
\end{pf}

\subsection{Variance and covariance of the integrated process}
\label{varcovar.subsec}

The goal of this section is to establish the asymptotic behavior of the
variance $V(Tt)=\operatorname{Var}[ \tilde{Y}_{Tt} ]$ as $T
\rightarrow\infty$ of the integrated process
%
\begin{equation}
\label{myy}
\tilde{Y}_t=\int_0 ^t Z(u)\,  du.
\end{equation}
Since $\{Z(u)\}$ is a Gaussian process with mean zero, so is $\{\tilde
{Y}_{t}\}$. Since $\{Z(u)\}$ is stationary, $\{\tilde{Y}_{t},  t\geq
0\}$ has stationary increments.

It is helpful to relate $V(t)$ exactly to the autocovariance $r(u)$.
This is the content of the following lemma. Notice in particular that
since $0<\bar{H}(s)<1$, $I_1$ is a finite constant for any fixed $X$
since it does not depend on $t$.

\begin{lemma}
\label{Vlemma}
For the variance $V(t)$, any constant $X \geq0$, and $t\geq X$
\[
V(t)=I_1+I_2+I_3,
\]
where $I_1 =2 \int_0^X\int_0^y r(x) \bar{H}_I(x) \, dx \, dy$,
$I_2 =2 c_X(t-X)$, $c_X=\int_0^X r(x) \bar{H}_I(x) \, dx$, and $I_3 =2
\int_X^t\int_X^y r(x) \bar{H}_I(x) \, dx \, dy$.
\end{lemma}

\begin{pf}
Since $Z(u)$ is a mean zero process it can be shown than
\[
V(t) =\mathbb{E}\biggl[ \int_0^tZ(u) \, du \int_0^t Z(v)\, dv \biggr]
=2\int_0^t \int_0^y r(x) \bar{H}_I(x) \, dx\, dy.
\]
Now separate the various regions of integration.
\end{pf}

The following theorem establishes the asymptotic results for the
variance $V(t)$ as $t \rightarrow\infty$ which is needed to prove the
main results.

\begin{lemma}
\label{Vttheorem}
Assume the autocovariance $r(u)$ satisfies the slowly varying function
condition in (\ref{regvarfuncr}). For the processes $\{Z(t), t\geq0\}
$ and $\{\tilde{Y}_t,  t\geq0\}$ defined in (\ref{myZ}) and (\ref
{myy}), respectively, and with
\[
\sigma^2= \frac{\sigma_{\mathit{lim}}^2(3-\alpha_\mathit{{min}})(2-\alpha
_\mathit{{min}})}{(\alpha_{\mathit{sess}}-1)(3-\alpha_\mathit{{min}}-\alpha_{\mathit{sess}})(4-\alpha
_\mathit{{min}}-\alpha_{\mathit{sess}})},
\]
Case \textup{3} ($1<\alpha_\mathit{{min}}<2$, $-1<2-\alpha_\mathit{{min}}-\alpha
_{\mathit{sess}} < 1$, $\alpha_{\mathit{sess}}>1$):
\[
V(t) \sim\sigma^2 t^{4-\alpha_\mathit{{min}}-\alpha_{\mathit{sess}}}L_V(t)L_r(t)
\qquad  \mbox{as } t \rightarrow\infty.
\]
Case \textup{4} ($1<\alpha_\mathit{{min}}<2$, $-2<2-\alpha_\mathit{{min}}-\alpha
_{\mathit{sess}}<-1$, $\alpha_{\mathit{sess}}<2$):
\begin{eqnarray*}
\lim_{T \rightarrow\infty} \frac{V(Tt)}{T L_V(T)L_r(T)} &=&\lim_{T
\rightarrow\infty} \frac{2ct}{L_V(T)L_r(T)} \qquad  \mbox{and}\\
V(t) &\sim&2ct \qquad  \mbox{as } t \rightarrow\infty \mbox{ where } c=\int
_0^\infty r(x)\bar{H}_I(x)\, dx.
\end{eqnarray*}
Boundary between cases 3 and 4 ($1<\alpha_\mathit{{min}}<2$, $2-\alpha
_\mathit{{min}}-\alpha_{\mathit{sess}}=-1$): Assume $c=\int_0^\infty r(x)\bar{H}_I(x)
\,dx < \infty$. Then
\[
\lim_{T \rightarrow\infty}\frac{V(Tt)}{T L_V(T)L_r(T)} = \infty
\quad  \mbox{and}\quad V(t) \sim2ct \qquad  \mbox{as } t \rightarrow\infty.
\]
\end{lemma}

\begin{pf}
By Corollary 1.4.2 of the Characterization theorem [\citet{bingham87}, pp. 17--18] there exists $X \geq0$ such that $L_V(x)$ and
$L_V(x)L_r(x)$ are locally bounded and locally integrable on $[X,\infty
)$. Since only the tail of $V(t)$ will be important, it is assumed $t >
X$. In particular, this means the integrals $I_1$ and $I_2$ are nonzero.

\textit{Case} 3:
By assumption $1<\alpha_{\mathit{sess}}<2$ and $\alpha_{\mathit{sess}}<3-\alpha
_\mathit{{min}}$. For $I_3$, let $J_1(y)$ be the inner integral
\[
J_1(y)=\int_X^y r(x) \bar{H}_I(x)\, dx
\]
so that $I_3=2\int_X ^t J_1(y)\,  dy.$
By (\ref{Hbartail}) and Karamata's theorem [\citet{bingham87}, p.~28],
$\mbox{as } t \rightarrow\infty$
\[
I_3\sim \frac{\sigma_{\mathit{lim}}^2 (3-\alpha_\mathit{{min}})(2-\alpha
_\mathit{{min}})}{(\alpha_{\mathit{sess}}-1)(3-\alpha_\mathit{{min}}-\alpha_{\mathit{sess}})(4-\alpha
_\mathit{{min}}-\alpha_{\mathit{sess}})}t^{4-\alpha_\mathit{{min}}-\alpha_{\mathit{sess}}}L_V(t)L_r(t).
\]
Since $4-\alpha_\mathit{{min}}-\alpha_{\mathit{sess}}>1$, $I_3$ is regularly varying
with higher order than $I_1$ and $I_2$ (which are constant and linear
in $t$, respectively), the result follows.

\textit{Case} 4:
By Karamata's theorem [\citet{bingham87}, p.~28]
%
\begin{eqnarray}
\label{int29}
I_2+I_3 &=& 2 \int_X^t \int_0^y r(x)\bar{H}_I(x) \,dx\, dy\nonumber\\[-8pt]\\[-8pt]
&=& 2 c (t-X) - 2 \int_X^t \int_y^\infty r(x)\bar{H}_I(x)\, dx\, dy,\nonumber
\end{eqnarray}
where $c=\int_0^\infty r(x)\bar{H}_I(x) \,dx < \infty$ since
$2-\alpha_\mathit{{min}}-\alpha_{\mathit{sess}}<-1$, $ \mbox{as } y \rightarrow
\infty$
\[
\int_y^\infty r(x)\bar{H}_I(x) \,dx \sim\frac{-\sigma
_{\mathit{lim}}^2(3-\alpha_\mathit{{min}})(2-\alpha_\mathit{{min}})}{2(\alpha
_{\mathit{sess}}-1)(3-\alpha_\mathit{{min}}-\alpha_{\mathit{sess}})}y^{3-\alpha_\mathit{{min}}-\alpha
_{\mathit{sess}}}L_V(y)L_r(y)
\]
and so for the integral in (\ref{int29})
\begin{eqnarray*}
&&\lim_{t\rightarrow\infty}\frac{ \int_X^t \int_y^\infty r(x)\bar
{H}_I(x)\, dx\, dy}{t^{4-\alpha_\mathit{{min}}-\alpha_{\mathit{sess}}}L_V(t)L_r(t)}  \\
&&\qquad=\frac{-\sigma_{\mathit{lim}}^2(3-\alpha_\mathit{{min}})(2-\alpha
_\mathit{{min}})}{2(\alpha_{\mathit{sess}}-1)(3-\alpha_\mathit{{min}}-\alpha_{\mathit{sess}}) (4-\alpha
_\mathit{{min}}-\alpha_{\mathit{sess}})}.
\end{eqnarray*}
Since $4-\alpha_\mathit{{min}}-\alpha_{\mathit{sess}}<1$ the linear term of (\ref
{int29}) is regularly varying with higher order, and $V(t) \sim2ct
 \mbox{ as } t \rightarrow\infty$
while
\[
\lim_{T\rightarrow\infty}\frac{V(Tt)}{Tt L_V(T) L_r(T)}= \lim
_{T\rightarrow\infty}\frac{2c}{L_V(T) L_r(T)}.
\]
The second limit requires more information about $L_V(T)$ and $L_r(T)$.
Both positive constants and the logarithm function are examples of
slowly varying functions, and they would give quite different limits.

\textit{Boundary case}:
By assumption, $2-\alpha_\mathit{{min}}-\alpha_{\mathit{sess}} =-1$ and $c=\int
_0^\infty r(x)\times \break \bar{H}_I(x) \, dx < \infty$. Since
\[
r(x)\bar{H}_I(x) \sim\frac{\sigma_{\mathit{lim}}^2}{2}\frac{1}{\alpha
_{\mathit{sess}}-1}(3-\alpha_\mathit{{min}})(2-\alpha_\mathit{{min}})x^{-1}L_V(x)L_r(x),
\]
we know $L_V(x)L_r(x) \rightarrow0$ as $x \rightarrow\infty$ since
otherwise would make the tail too heavy for integrability.

Let $J_2(y)=\int_y^\infty r(x)\bar{H}_I(x) \, dx$ and $c_X=\int_0^X
r(x)\bar{H}_I(x)\, dx$. Then
\[
V(t)= I_1 + 2 \int_X^t\bigl(c-J_2(y)\bigr)\,  dy,
\]
and so the asymptotic properties of $V(t)$ depend, in part, on those of
$J_2(y)$ as $y \rightarrow\infty$. For $J_2(y)$, $\lim_{y\rightarrow
\infty}J_2(y) = 0$, by [\citet{bingham87}, Proposition 1.5.9(b)]
$J_2(y)$ is slowly varying with
\[
\lim_{y\rightarrow\infty}\frac{J_2(y)}{L_V(y)L_r(y)}= \infty,
\]
and by Karamata's theorem [\citet{bingham87}, p. 28] it follows that
\begin{eqnarray*}
\lim_{t \rightarrow\infty}\frac{\int_X^t J_2(y)\, dy}{t
J_2(t)}&=&1,\qquad   \lim_{t\rightarrow\infty}\frac{\int_X^t J_2(y) \,dy}{t
L_V(t)L_r(t)}=\infty  \quad  \mbox{and}\\
 \lim_{t \rightarrow\infty
}\frac{\int_X^t J_2(y) \,dy}{t}&=&0.
\end{eqnarray*}
The results follow immediately.
\end{pf}

\subsection{Proof of main results}
\label{mainproof.subsec}

\begin{pf*}{Proof of Theorem \ref{fbmmainthm1}}
By Theorem \ref{incrementalprocess},
\[
\mathcaligr{L}\lim_{n \rightarrow\infty} \frac{\sum_{i=-\infty
}^{\infty}X_{n,i}(t)}{\sqrt{\lambda_n}}= Z(t)
\]
is Gaussian, mean zero and has covariance given by $\operatorname{Cov}[Z(s),Z(t)]=r(t-s) \int_{t-s} ^{\infty} \bar{H}(z)\,  dz,$ $s<t$.
Now, for any $t\geq0$ and $T>0$ let $Y(t)$ be defined by
\[
Y(t)= \frac{\tilde{Y}_{t}}{ \sqrt{L_V(T)L_r(T)}T^{H}}.
\]
Then for the characteristic function
\begin{eqnarray*}
\mathbb{E}[\exp(i s Y(Tt))] &=& \exp\biggl(-\frac{V(Tt)}{2
L_V(T)L_r(T)T^{2H}} s^2\biggr)\\
&\rightarrow&\exp\biggl(- \frac{\sigma_{\mathit{lim}}^2(2H-1+\alpha
_{\mathit{sess}})}{2H(2H-1)} s^2\biggr) \qquad  \mbox{as } T \rightarrow\infty
\end{eqnarray*}
since
\[
\lim_{T\rightarrow\infty}\frac{V(Tt)}{ L_V(T)L_r(T)T^{2H}}= \frac
{\sigma_{\mathit{lim}}^2(2H-1+\alpha_{\mathit{sess}})}{2H(2H-1)}  t^{2H}
\]
by Lemma \ref{Vttheorem}.
Thus, for fixed $t\geq0$
\[
Y(Tt) \stackrel{D}{\rightarrow} Y_1 \qquad  \mbox{as }
T\rightarrow\infty,
\]
where
\[
Y_1 \sim\mathrm{N}\biggl(0,\frac{\sigma_{\mathit{lim}}^2(2H-1+\alpha
_{\mathit{sess}})}{2H(2H-1)}t^{2H}\biggr).
\]

For the covariances, by Lemma \ref{Vttheorem} and stationary increments
\begin{eqnarray*}
&&\operatorname{Cov}[Y(Ts),Y(Tt)] \\
&&\qquad = \frac{1}{2} \frac{(V(Ts) + V(Tt)-V(T(t-s))
)}{L_V(T)  L_r(T)  T^{4-\alpha_\mathit{{min}}-\alpha_{\mathit{sess}}} }\\
&&\qquad \rightarrow\sigma^2 \bigl(s^{4-\alpha_\mathit{{min}}-\alpha_{\mathit{sess}}} +
t^{4-\alpha_\mathit{{min}}-\alpha_{\mathit{sess}}} - (t-s)^{4-\alpha_\mathit{{min}}-\alpha
_{\mathit{sess}}} \bigr) \qquad  \mbox{as } T \rightarrow\infty\\
&&\qquad = \sigma^2 \bigl(s^{2H} + t^{2H}-(t-s)^{2H} \bigr).
\end{eqnarray*}
Thus, by the definition of fractional Brownian motion [\citet{st94}] the
left and right-hand sides of (\ref{fbmmain1eqn}) have the same
finite-dimensional distributions.
\end{pf*}

\begin{pf*}{Proof of Theorem \ref{fbmmainthm2}}
The proof is similar to that for Theorem \ref{fbmmainthm1} but shows, for
\[
Y(t) = \frac{\tilde{Y}_t}{T^{1/2}}
\]
and fixed $t\geq0$
\[
Y(Tt) \stackrel{D}{\rightarrow} Y_2 \qquad  \mbox{as } T\rightarrow\infty
\mbox{ where }
Y_2 \sim\mathrm{N}(0,2ct)
\]
and for $s<t$
\[
\operatorname{Cov}[Y(Ts),Y(Tt)] = \frac{1}{2} \frac{(V(Ts) +
V(Tt)-V(T(t-s)) )}{T}
\rightarrow2cs \qquad  \mbox{as } T \rightarrow\infty
\]
so the left and right-hand sides of (\ref{fbmmainthm2eqn}) have the
same finite-dimensional distributions.
\end{pf*}

\section{Weak convergence}
\label{weakconv:sec}

In this section, weak convergence is established for both the first
limit ($n \rightarrow\infty$) providing the incremental process, and
the second limit ($T \rightarrow\infty$) providing the limiting process.

\subsection{Weak convergence for the first limit}

\begin{theorem}
The convergence in Theorem \ref{incrementalprocess} can be
strengthened to weak convergence in the space $D[0,\infty)$ equipped
with the $J_1$ topology, and the limiting Gaussian process is almost
surely continuous.
\end{theorem}

\begin{pf}
Convergence of the finite-dimensional distributions was established in
Theorem \ref{incrementalprocess} so it remains to prove tightness. For
any $M>0$, take $0 \leq t_1 \leq t_2 \leq t_3 \leq M$, $t_3-t_1<1$ and define
%
\begin{eqnarray}
\label{weakconv1}
U_n(t)&=&\frac{\sum_{i=1}^{N_n(t)} X_{n,i}(t)}{\sqrt{\lambda_n}}
\quad  \mbox{and}\nonumber\\
[-8pt]\\[-8pt]
 \Delta&=&\bigl(U_n(t_2)-U_n(t_1)\bigr)^2
\bigl(U_n(t_3)-U_n(t_2)\bigr)^2.\nonumber
\end{eqnarray}
To prove tightness in $D[0,\infty)$ it will be shown that there exist
constants $K_1>0$ and $K_2>0$ (possibly depending on $M$) such that
%
\begin{eqnarray}
\label{weak1}
E[\Delta] &\leq& K_1(t_3-t_1)^2,\qquad   0<t_3-t_1 < 1\quad  \mbox{and}\\
\label{weak2}
E[U_n(t_2)-U_n(t_1)] &\leq& K_2(t_3-t_1),\qquad   0<t_3-t_1 < 1.
\end{eqnarray}
Then the result will follow by a result due to \citeauthor{whitt02} (\citeyear{whitt02}, pp.~226--227).

First we prove the bound of (\ref{weak1}). Establishing the bound of
(\ref{weak2}) will require little additional work. Let $
(X_{n,k}^{(i,j)}(t_{i+1}),\ldots,X_{n,k}^{(i,j)}(t_{j}) )$ and
$N_{n,i,j}(t)$ be defined as in the proof of Theorem \ref
{incrementalprocess} (so mutually independent in $i$,$j$, and $k$). For
simplicity define $N_{n,i,j}=N_{n,i,j}(t_3)$. Define
\[
U_{i,j,k}=\frac{1}{\sqrt{\lambda_n}}\sum
_{u=1}^{N_{n,i,j}}X_{n,u}^{(i,j)}(t_k)
\]
and
\[
T_{i,j,k}=\frac{1}{\sqrt{\lambda_n}}\sum_{u=1}^{N_{n,i,j}}
\bigl(X_{n,u}^{(i,j)}(t_k) -X_{n,u}^{(i,j)}(t_{k-1})\bigr)=U_{i,j,k}-U_{i,k,(k-1)}.
\]
Since $E[X_{n,u}^{(i,j)}(t) ]=0$, by conditioning arguments it follows
that $E[U_{i,j,k}]=E[T_{i,j,k}]=0$. In what follows below both
$U_{i,j,k}$ and $T_{i,j,k}$ will contribute to $E[\Delta]$ through the
nonzero higher moments of the summands. In addition, $T_{i,j,k}$ will
contribute through the temporal correlation within each summand.

It can be shown that
\begin{eqnarray*}
U_n(t_2)-U_n(t_1) &=&T_{0,2,2}+T_{0,3,2}+U_{1,2,2}+U_{1,3,2}-U_{0,1,1}
\qquad  \mbox{and}\\
U_n(t_3)-U_n(t_2) &=&T_{0,3,3}-T_{1,3,3}+U_{2,3,3}-U_{0,2,2}-U_{1,2,2}.
\end{eqnarray*}
Then the expansion of $\Delta$ in (\ref{weakconv1}) has $25^2$ terms
before simplifications. Those terms whose factors differ in either the
first or second indices ($i$ or $j$) are the product of independent
terms. Since an individual term has expected value zero, $E[\Delta]$
gets contributions from only 31 nonzero terms. That is,
%
\begin{eqnarray}
\label{weakYexpansion}
E[\Delta]&=&  E[
(U_{1,2,2}^2+U_{1,3,2}^2+T_{0,3,2}^2+U_{0,1,1}^2+T_{0,2,2}^2)\nonumber\\
&&\hspace*{9pt}{}\times(U_{1,2,2}^2+U_{0,2,2}^2+U_{2,3,3}^2+ T_{1,3,3}^2+T_{0,3,3}^2
)]\nonumber\\
&& {}+4E[-T_{0,2,2}U_{1,3,2}T_{1,3,3}U_{0,2,2} -
T_{0,3,2}U_{1,2,2}T_{0,3,3}U_{1,2,2} \\
&&\hspace*{29pt}{}-U_{1,2,2}U_{1,3,2}T_{1,3,3}U_{1,2,2}+T_{0,2,2}U_{1,2,2}U_{0,2,2}U_{1,2,2}\nonumber\\
&&{}\hspace*{29pt}-T_{0,2,2}T_{0,3,2}T_{0,3,3}U_{0,2,2}
+T_{0,3,2}U_{1,3,2}T_{0,3,3}T_{1,3,3}].\nonumber
\end{eqnarray}
The idea now will be to bound both $\sigma_W^2-r(u)$ and
$p_{i,j}(t_3)$ (recall $E[N_{n,i,j}]=\lambda_n p_{i,j}(t_3)$) suitably.

From now $K$ will be an arbitrary constant whose exact value will vary
from one calculation to the next. By the assumed Lipschitz continuity
of $H(x)$ we have
\begin{eqnarray*}
p_{1,2}(t_3)&=&\int_{t_1}^{t_2}H(t_3-s)-H(t_2-s) \,ds\\
&\leq&\int
_{t_1}^{t_2} K(t_3-t_2) \, ds
\leq K(t_3-t_1)^2.
\end{eqnarray*}
Since $0 \leq H(x) \leq1$, $p_{0,1}(t_3)$, $p_{0,2}(t_3)$,
$p_{1,3}(t_3)$, $p_{2,3}(t_3)$ are bounded above by $(t_3-t_1)$. Thus
we have
\begin{eqnarray*}
E[U_{1,2,2}^2] &=&\frac{1}{\lambda_n}E\Biggl[\sum
_{u=1}^{N_{n,1,2}}X_{n,u}^{(1,2)}(t_2)\sum
_{v=1}^{N_{n,1,2}}X_{n,v}^{(1,2)}(t_2) \Biggr] \nonumber\\
&=&\frac{1}{\lambda_n}E[N_{n,1,2}\sigma_W^2] \nonumber\\
&\leq& K(t_3-t_1)^2
\end{eqnarray*}
and for $(i,j) \neq(1,2)$
%
\begin{equation}
\label{boundUk}
E[U_{i,j,k}^2]\leq K p_{i,j}(t_3)\leq K(t_3-t_1).
\end{equation}
Since $E[N_{n,i,j}^2-N_{n,i,j}]=\lambda_n^2 p_{i,j}^2$ it follows that
\[
E[U_{1,2,k}^4]=\frac{1}{\lambda
^2}E\bigl[N_{n,1,2}E[X_{n,i}^4(t_k)]+3(N_{n,1,2}^2-N_{n,1,2})\sigma_W^4 \bigr]
\leq K(t_3-t_1)^2.
\]

To understand contributions from $T_{i,j,k}$ notice that we can write
%
\begin{equation}
\label{weakcor}
\sigma^2_W-r(t)=r(0)-r(t)=\frac{\mu_1}{\mu_1+\mu_2}\bigl(1-\pi_{11}(t)\bigr),
\end{equation}
where $\pi_{11}(t)=P(W(t)=1|W(0)=1)$, which is continuous in $t$ since
the on and off distributions are assumed continuous. Now by \citet
{taqqu97}, equation (13),
\[
1-\pi_{11}(t)=\frac{1}{\mu_1}\int_0^t\bar{F}_1(u)\, du+\int_0^t
\bar{F}_1(t-u)h_{12}(u)\, du,
\]
where $H_{12}$ is the renewal function for the interrenewal
distribution $F_1*F_2$ (i.e., $H_{12}=\sum_{k=1}^\infty(F_1*F_2)^k$)
with the density $h_{12}(u)$. Using the Fundamental Theorem of Calculus
for the first term, and Laplace transform arguments for the second term,
\[
\frac{d}{dt}\bigl(1-\pi_{11}(t)\bigr)=\frac{1}{\mu_1}\bar{F}_1(t)-\int_0^t
\bar{f}_1(t-u)h_{12}(u) \,du.
\]
On the right side, the first term is bounded between 0 and $1/\mu_1$.
The second term is also bounded [\citet{bhat}, p. 167]. Therefore, since
$(1-\pi_{11}(t))$ is continuous with a bounded density, it is
Lipschitz continuous. Using (\ref{weakcor}) it follows that
\[
|\sigma_W^2-r(u) | \leq Ku.
\]
It can also be shown that
\begin{eqnarray*}
E[U_{i,j,k}U_{i,j,(k-1)} ]&=&\frac{1}{\lambda_n}E\Biggl[\sum
_{u=1}^{N_{n,i,j}}X_{n,u}^{(i,j)}(t_k)\sum
_{v=1}^{N_{n,i,j}}X_{n,v}^{(i,j)}\bigl(t_{(k-1)}\bigr) \Biggr] \\
&=&p_{i,j}(t_3) r\bigl(t_k-t_{(k-1)}\bigr)
\end{eqnarray*}
and by expanding the square
%
\begin{eqnarray}
\label{tightT}
E[T_{i,j,k}^2]&=&E\bigl[\bigl(U_{i,j,k}-U_{i,j,(k-1)}\bigr)^2\bigr] \nonumber\\
&\leq&2p_{i,j}(t_3)\big|\sigma_W^2-r\bigl(t_k-t_{(k-1)}\bigr)\big| \\
&\leq& K(t_3-t_1)^2.\nonumber
\end{eqnarray}
Notice that since $|\sigma_W^2-r(u)|\leq2 \sigma_W^2$ the looser
bound $E[T_{i,j,k}^2]\leq K p_{i,j}(t_3)$ is also available.

Exploiting the idea that since $\{W(t)\}$ is a 0--1 process we know for
the joint moments $E[W^r(t_i)W^s(t_j)]=E[W(t_i)W(t_j)]$ for $r$, $s \in
\mathbb{Z}^+$ it can also be shown that
\[
E[U_{1,3,2}U_{1,3,3}^3 ]=\frac{p_{1,3}(t_3)}{\lambda_n}(1-3\mu+3\mu
^2)r(t_3-t_2)+3p^2_{1,3}(t_3)\sigma_W^2r(t_3-t_2)
\]
and
\begin{eqnarray*}
E[U_{1,3,2}^2 U_{1,3,3}^2] &= & \frac
{p_{1,3}(t_3)}{\lambda_n} (1-2\mu^2) r(t_3-t_2) +\frac
{p_{1,3}(t_3)}{\lambda_n}\sigma_w^4 + p^2_{1,3}(t_3)\sigma_w^4 \\
&&{}+2p^2_{1,3}(t_3)[r(t_3-t_2)]^2.
\end{eqnarray*}
Fourth moments of $T_{i,j,k}$ can now be found. For example, we have that
\begin{eqnarray*}
E[T_{1,3,3}^4]&=&E[U_{1,3,3}^4-4U_{1,3,3}^3
U_{1,3,2}+6U_{1,3,3}^2U_{1,3,2}^2-4U_{1,3,3}U_{1,3,2}+U_{1,3,2}^4]\\
&=&\frac{2p_{1,3}(t_3)}{\lambda_n}\bigl(\sigma
_W^2-r(t_3-t_2)\bigr)+12p_{1,3}(t_3)\sigma_W^2\bigl(\sigma_W^2-r(t_3-t_2)\bigr)\\
&&{}+12p^2_{1,3}\bigl(\sigma_W^4-r^2(t_3-t_2)\bigr)
\\
&\leq&\frac{2p_{1,3}(t_3)}{\min_n \{\lambda_n\}
}K(t_3-t_2)\\
&&{}+12[p_{1,3}(t_3)\sigma_W^2K(t_3-t_2)+p^2_{1,3}(t_3)(2\sigma_W^2)K(t_3-t_2)]
\\
&\leq& K(t_3-t_1)^2.
\end{eqnarray*}

The inequalities above provide a means to bound the first 25 terms in
the expansion of (\ref{weakYexpansion}). For the remaining six terms
we use the Cauchy--Schwarz inequality and independence. For example, it
can be shown that
\begin{eqnarray*}
E[T_{0,2,2}U_{1,3,2}T_{1,3,3}U_{0,2,2}] &\leq&(E[T_{0,2,2}^2
U_{1,3,2}^2])^{1/2} (E[T_{1,2,3}^2 U_{0,2,2}^2]
)^{1/2}\\
&=&(E[T_{0,2,2}^2] E[U_{1,3,2}^2] E[T_{1,3,3}^2]
E[U_{0,2,2}^2])^{1/2}\\
&\leq& K(t_3-t_1)^2.
\end{eqnarray*}
Now every term in the expansion of (\ref{weakYexpansion}) is bounded
by $K(t_3-t_1)^2$ for some constant $K>0$ and (\ref{weak1}) is shown.

To establish (\ref{weak2}) note that for $t_3-t_1<1$
\begin{eqnarray*}
E\bigl[\bigl(U_n(t_3)-U_n(t_1)\bigr)^2\bigr]&=&E[(T_{0,3,3} + U_{1,3,3} +
U_{2,3,3} - U_{0,1,1}- U_{0,2,1})^2] \\
&\leq& KE[(T_{0,3,3}^2 + U_{1,3,3}^2 + U_{2,3,3}^2 + U_{0,1,1}^2 +
U_{0,2,1}^2)] \\
&\leq& K[(t_3-t_1)^2 + (t_3-t_1)] \\
&\leq& K_2 (t_3-t_1),
\end{eqnarray*}
where the second line follows from independence for distinct indices
and the third line follows from (\ref{boundUk}) and (\ref{tightT}).
Thus, weak convergence in $D[0,\infty)$ and limiting almost surely
continuous paths are established [\citet{whitt02}, pp. 126--127].
\end{pf}

\subsection{Weak convergence for the second limit}

\begin{theorem}
The convergence (as $T \rightarrow\infty$) in Theorems \ref
{fbmmainthm1} and \ref{fbmmainthm2} can be strengthened to weak
convergence in $C[0,\infty)$ with the uniform topology.
\end{theorem}

\begin{pf}
Recall the integrated process
\[
\tilde{Y}_t=\int_0^t Z(u)\, du
\]
defined in Section \ref{varcovar.subsec} which has stationary
increments. Since convergence of the finite-dimensional distributions
is established above (Theorems \ref{fbmmainthm1} and \ref
{fbmmainthm2}) it remains to prove tightness, which we do using a
moment condition [\citet{billingsley68}, p. 95].

\textit{Case} 3:
It will be shown that for $0<u<1$ and $T$ sufficiently large
\[
E\biggl[\frac{(\tilde{Y}_{T(s+u)}-\tilde
{Y}_{Ts})^2}{T^{2H}L_V(t)L_r(t)} \biggr]=\frac
{V(Tu)}{T^{2H}L_V(t)L_r(t)}\leq Ku^{1+\delta}
\]
for some constant $K>0$ and some constant $\delta>0$.

Since $V(t) \sim\sigma^2 t^{2H}L_V(t)L_r(t)$ we have
\[
\frac{V(Tu)}{T^{2H}L_V(t)L_r(t)} \leq\frac
{K_1(Tu)^{2H}L_V(Tu)L_r(Tu)}{T^{2H}L_V(T)L_r(T)}.
\]
Proceeding as in [\citet{mikosch2002}] (or using the Potter bounds
[\citet{bingham87}, p. 25]), there exists $t_0$ such that for all $u
\leq1$ and $Tu \geq t_0$
\[
\frac{L_V(Tu)L_r(Tu)}{L_V(T)L_r(T)} \leq K_2 u^{-\delta}
\]
with small $\delta>0$ chosen so that $2H-2\delta>1$. Then
\[
\frac{V(Tu)}{T^{2H}L_V(T)L_r(T)}\leq K_1 K_2 u^{2H-\delta}.
\]
For $0<Tu<t_0$, since $V(t) \sim\sigma^2 t^{2H}L_V(t)L_r(t)$ for $T$
sufficiently large we have
\begin{eqnarray*}
\frac{V(Tu)}{T^{2H}L_V(T)L_r(T)} &\leq&\frac{2
(Tu)^2}{T^{2H}L_V(T)L_r(T)}\leq\frac{2 (Tu)^{1+\delta
}}{T^{2H}L_V(T)L_r(T)}t_0^{1-\delta}u^{1+\delta}\\
& \leq& K_3 u^{1+\delta
}.
\end{eqnarray*}
Since $2H-\delta<1+\delta$ we have
\[
\frac{V(Tu)}{T^{2H}L_V(T)L_r(T)}\leq Ku^{1+\delta}
\]
for $T$ large enough.

\textit{Boundary case and case} 4:
Since $\{ \tilde{Y}_t,  t \geq0 \}$ has stationary increments it is
enough to show for that for $u>0$ and some constant $K$
\[
E\biggl[ \biggl(\frac{\tilde{Y}_{Tu}}{T^{1/2}} \biggr)^4 \biggr]
\leq K u^2.
\]
Now, we know $\tilde{Y_t} \sim N(0,V(t))$ and so $E[\tilde
{Y}_t^4]=3V^2(t)$, while
\[
V(t) = 2ct -2\int_0^t \int_y ^\infty r(x) \bar{H}_I \, dx \,  dy \leq
2ct.
\]
Thus for all $T>0$ and $u>0$ and some constant $K$
\[
E\biggl[ \biggl(\frac{\tilde{Y}_{Tu}}{T^{1/2}} \biggr)^4 \biggr] =
3 \frac{V^2(Tu)}{T^2} \leq K u^2,
\]
tightness follows, and weak convergence is established [\citet{billingsley68}, p.~95].
\end{pf}

\section{Simulation results}
\label{simresults:sec}

In this section the limit results obtained above are demonstrated
through simulation. The starting point for each simulation is the
``session'' which arrives according to a Poisson process with intensity
$\lambda$. In these simulations each \mbox{session} has a Pareto distributed
lifetime with mean $\mu_{\mathit{sess}}$ and characteristic exponent $\alpha
_{\mathit{sess}}$, $1<\alpha_{\mathit{sess}}<2$. To simplify the implementation, session
lengths were truncated to integers. This step produced a list of start
and end times for the ``sessions.'' In a second step, these start and
end times were converted to a timeseries $\{B(k), k=1,\ldots,25\times
10^6\}$, the number of sessions alive at time $k$. (For these
simulations, the first $50\times10^6$ observations were discarded to
achieve stationarity and the next $25 \times10^6$ observations were
retained.) This corresponds to the number of busy servers in the
infinite source Poisson model (i.e., the busy server process of an
$M/G/\infty$ queue).

On--off activity within the sessions was simulated in a third step.
Starting with the same session start and end times, session activity
within each session was simulated as a stationary on--off process with
parameters $\alpha_1$, $\mu_1$, $\alpha_2$ and $\mu_2$. This
produced a dataset $\{A(k), k=1,\ldots,25\times10^6\}$ which
represents the number of active (i.e., alive and on) sources for each
$k$. Finally, a dataset $\{C(k), k=1,\ldots,25 \times10^6\}$ was
produced where
\[
C(k)=\frac{\mu_1}{\mu_1+\mu_2}B(k)-A(k)
\]
and corresponds to the increments of (\ref{randomcentering}) above.

As shown in Table \ref{simtable}, several series, each corresponding
to different parameter values, were simulated. For the datasets in
series 1 through series 3, five simulations were performed in each
series. For series 4, 10 simulations were performed. For all
simulations, the parameters for the on--off activity were unchanged,
with $\alpha_1=\alpha_2=1.4$ and $\mu_1=\mu_2=100$. The parameters
for the sessions, namely $\alpha_{\mathit{sess}}$, $\mu_{\mathit{sess}}$, and $\lambda
$ are shown in Table \ref{simtable}. The values of $\alpha_1$,
$\alpha_2$ and $\alpha_{\mathit{sess}}$ were chosen mainly for illustrative
purposes. Hurst parameters can be calculated for the infinite source
Poisson model ($H_{\mathit{ISP}}=(3-\alpha_{\mathit{sess}})/2$), the stationary on--off
model [\citet{taqqu97}] ($H_{\mathit{on}/\mathit{off}}=(3-\min(\alpha_1,\alpha_2))/2$),
and the Hurst parameter $H_{\mathit{burst}/\mathit{on}/\mathit{off}}$ described above in Theorems
\ref{fbmmainthm1} or \ref{fbmmainthm2}. The theoretical values are
all shown in the table.

\begin{table}[b]
\caption{Parameter values for the simulations}\label{simtable}
\begin{tabular*}{\textwidth}{@{\extracolsep{\fill}}l c c c c c c c c@{}}
\hline
& $\alpha_{\mathit{sess}}$ & $\mu_{\mathit{sess}}$ & $\lambda$ & $\mu_1$, $\mu_2$ &
$\alpha_1$, $\alpha_2$ & $H_{\mathit{burst}/\mathit{on}/\mathit{off}}$ & $H_{\mathit{ISP}}$ &
$H_{\mathit{on}/\mathit{off}}$ \\
\hline
Series 1 & 1.2 & 12000 & 1 & 100 & 1.4 & 0.7 & 0.9 & 0.8 \\
Series 2 & 1.2 & \phantom{00}120 & 5 & 100 & 1.4 & 0.7 & 0.9 & 0.8 \\
Series 3 & 1.2 & \phantom{0}1200 & 1 & 100 & 1.4 & 0.7 & 0.9 & 0.8 \\
Series 4 & 1.8 & \phantom{0}1200 & 1 & 100 & 1.4 & 0.5 & 0.6 & 0.8 \\
\hline
\end{tabular*}
\end{table}

Hurst parameter estimates for the simulated data were made using the
``logscale diagram'' [\citet{abry98}]. This technique divides the data
into blocks of size $2^j$ and calculates the variance of the $n$th
order Daubechies wavelet coefficients for each value of $j$. Values of
$j$ are sometimes called ``octaves.'' Collectively, the points show how
the data scales with increasing block length. $H$ is estimated using
weighted linear regression to find the slope of the graph over an
interval of octaves chosen by the user. Besides being quick, this
technique provides approximate 95\% confidence intervals.

\begin{figure}

\includegraphics{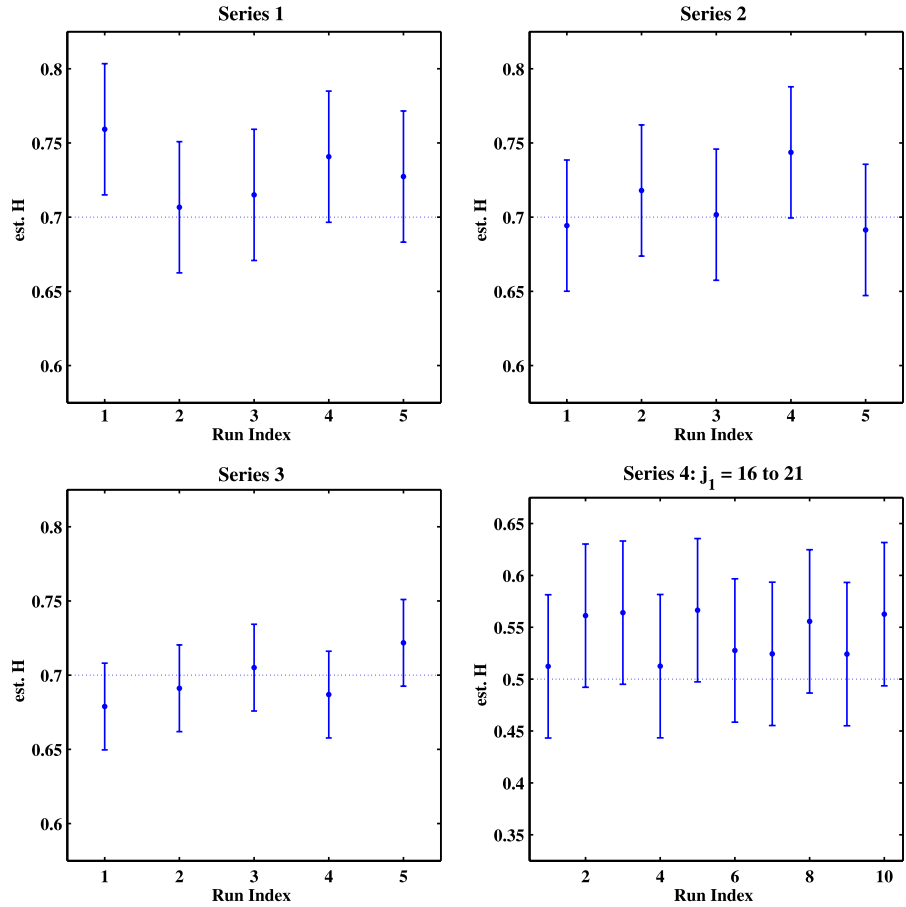}

\caption{Estimated Hurst parameters with
approximate 95\% confidence intervals for the simulated data.}\label{simfigures}
\end{figure}

Figure \ref{simfigures} shows the results from series 1--4. For each
run, the central dot is the value of the estimated Hurst parameter,
while the high and low bars show the ends of the approximate 95\%
confidence inteval. For series 1 (top left), the estimate is peformed
over the interval $j=15,\ldots,21$. Since the mean session length is
12,000, and $\log_2 12{,}000 \sim13.6$, $j=15$ provides a starting
octave beyond any artifacts associated with the mean session length.
For series 2 (top right) and series~3 (bottom left), with shorter mean
session lengths, $j=14,\ldots,21$ was used for the regression. In all
cases the estimated Hurst parameter is within 0.06 of the value
described in Theorem \ref{fbmmainthm1}, and in all but one case the
theoretical value is within the confidence interval.

\begin{figure}

\includegraphics{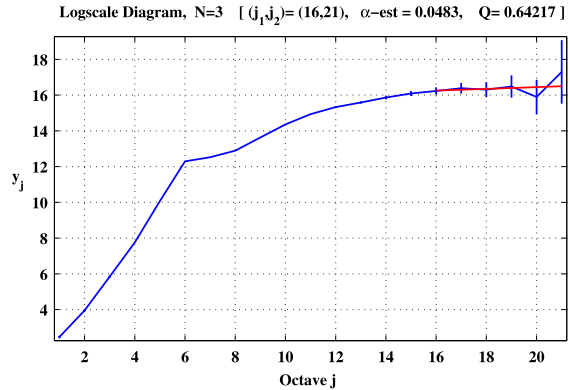}

\caption{Logscale diagram for run 8 of series 4
(right). An ultimately horizontal line corresponds to $H=0.5$, and is
predicted by Theorem \textup{\protect\ref{fbmmainthm2}}.}\label{series4figures}
\end{figure}

Figure \ref{simfigures} (bottom right) shows the Hurst parameter
estimates for the series 4 data using octaves 16 to 21. The horizontal
line corresponds to a value of 0.5, which is the Hurst parameter
predicted by Theorem \ref{fbmmainthm2}. The estimated values are
somewhat higher than that predicted by theory, although the confidence
intervals do cover the theoretical value. The logscale diagram in
Figure \ref{series4figures} shows the difficulties in estimating the
Hurst parameters with this data. One hopes to find an interval of
linear increase on which to obtain the Hurst parameter estimate, but
where does such an interval start? A horizontal line would correspond
to a Hurst parameter of~0.5. A reasonable start is octave 16 for which
the estimate would be 0.524. But, it would be equally reasonable to
start at a larger octave, for which the estimate would be even closer
to 0.5. All the data for series 4 is similar, and so estimates using
octaves 16 to 21 are shown here. Moreover, since the trend in each
logscale diagram is to flatness, the point estimates are likely to be
consistently biased too large. Even so, the results from these
simulations appear to agree with the value predicted in Theorem \ref
{fbmmainthm2}. In all cases the confidence interval covers the
predicted value.

\section{Conclusions}
\label{conc:sec}

In this paper, the infinite source Poisson model for bursts has been
combined with the on--off model for in-burst activity. Burst and on--off
durations are assumed heavy tailed with infinite variance and finite
mean. By using convergence results for random centering of random sums,
weak convergence of the workload to fractional Brownian motion is
shown. The degree of long-range dependence is shown to depend on the
tail indices of both the on--off durations and the lifetimes
distributions. Moreover, the results can be separated into cases
depending on those tail indices. In one case where all distributions
are heavy tailed it is shown that the limiting result is Brownian motion.

The method of proof here (and in \citet{rolls2003}) uses iterated
limits for first the arrival rate, and then the time scale, as in
\citet{taqqu97} and \citet{brichet00}. As in those papers, one could
study convergence with the limits taken in reverse order. This is left
for future work. Moreover, the iterated limits are not entirely
satisfactory, and a single limit as in \citet{mikosch2002} and \citet
{gagailas03} would be preferable. Unfortunately, the change in details
within the proof would make it essentially a new proof, and is left for
future work.

Random centering in traffic models appears to be a new idea. The
application is not necessarily targeted at current TCP/IP data networks
but simply any network where the start and end of some kind of
``sessions'' are signaled. It would be interesting to see if in a
network element this can be practically used to advantage, possibly
leveraging the ``nicer'' queueing behavior of Brownian motion workloads.


\section*{Acknowledgments}

Portions of this work were supported by a Queen's University Thesis
Completion Bursary or a UNC-Wilmington Summer Research Initiative
grant. I am also indebted to my thesis supervisor, Dr. Glen K.
Takahara, for helpful conversations when this work was in its infancy.
Finally I would like to thank the anonymous referees for their comments
which improved this paper.
%

\printaddresses

\end{document}